\documentstyle[12pt]{article}
\begin{document}

\newcommand{\singlespace}
{

\renewcommand{\baselinestretch}{1}
\large\normalsize}

\newcommand{\doublespace}{

\renewcommand{\baselinestretch}{1.2}

\large\normalsize}

\renewcommand{\theequation}{\thesection.\arabic{equation}}
\input amssym.def
\input amssym
\setcounter{equation}{0}
\def \ten#1{_{{}_{\scriptstyle#1}}}

\def \g{\frak g}

\def \h{\frak h}

\def \1{\bf 1}

\def \Z{\Bbb Z}

\def \C{\Bbb C}

\def \R{\Bbb R}

\def \Q{\Bbb Q}

\def \N{\Bbb N}

\def \l{\lambda}

\def \V{V^{\natural}}

\def \wt{{\rm wt}}

\def \tr{{\rm tr}}

\def \Res{{\rm Res}}

\def \End{{\rm End}}

\def \Aut{{\rm Aut}}

\def \mod{{\rm mod}}

\def \irr{{\rm irr}}

\def \Hom{{\rm Hom}}

\def \im{{\rm im}}

\def \<{\langle}

\def \>{\rangle}

\def \w{\omega}

\def \o{\omega}

\def \t{\tau }

\def \ch{{\rm ch}}

\def \a{\alpha }

\def \b{\beta}

\def \e{\epsilon }

\def \la{\lambda }

\def \om{\omega }

\def \O{\Omega}

\def \qed{\mbox{$\square$}}

\def \pf{\noindent {\bf Proof: \,}}

\def \voa{vertex operator algebra\ }

\def \voas{vertex operator algebras\ }

\def \p{\partial}

\def \1{{\bf 1}}

\singlespace

%\doublespace

\newtheorem{thmm}{Theorem}
\newtheorem{th}{Theorem}[section]

\newtheorem{prop}[th]{Proposition}

\newtheorem{coro}[th]{Corollary}

\newtheorem{lem}[th]{Lemma}

\newtheorem{rem}[th]{Remark}

\newtheorem{de}[th]{Definition}

\newtheorem{con}[th]{Conjecture}

\newtheorem{ex}[th]{Example}
\begin{center}
{\Large{\bf  Integrability of $C_2$-Cofinite Vertex Operator Algebras}}
\markboth{Integrability of $C_2$-Cofinite VOAs}
{Integrability of $C_2$-Cofinite VOAs}
\vspace{0.5cm}

Chongying Dong\footnote{Supported by NSF grants, China NSF grant 
10328102 and faculty research funds
granted by the University of
California at
Santa Cruz.}
and

Geoffrey Mason\footnote{Supported by a NSF grant and faculty
research funds granted by the University of
California at
Santa
Cruz.}
\\

Department of Mathematics, University of
California,
Santa Cruz, CA 95064
\end{center}

\hspace{1.5
cm}

\begin{abstract}
\noindent
The following \emph{integrability} theorem for
vertex operator algebras $V$ satisfying some finiteness conditions
($C_2$-cofinite
and  CFT-type) is proved:  the
vertex operator subalgebra generated
by a simple Lie subalgebra $\g$ of the
weight one subspace $V_1$ is isomorphic to
the irreducible highest weight $\hat{\g}$-module $L(k, 0)$ for a 
positive integer $k$, and $V$ is an integrable
  $\hat{\g}$-module. The case in which $\g$
is replaced by an abelian Lie subalgebra is also considered, and several
consequences of integrability are discussed.
  2000MSC:17B69
\end{abstract}

\section{Introduction}
\setcounter{equation}{0}

Representations of affine Kac-Moody algebras
by vertex operators not
only partially motivated the idea of a vertex (operator) algebra,
but
also provide an important class of vertex operator algebras (cf. [FZ]). In this
paper
we are concerned with affine Kac-Moody algebras contained in an
abstract vertex operator
algebra.

  The vertex operator algebras
$V=(V,Y,\1,\o)$ that we study are \emph{simple} and of \emph{strong
CFT-type}. Following [DLMM], [DM1], we say that $V$ is of \emph{CFT-type} in
case the $L(0)$-grading on $V$ has the shape
\begin{eqnarray*}
V= \bigoplus_{n \geq 0}V_n = \C\1 \oplus V_1 \oplus ...
\end{eqnarray*}
It has \emph{strong} CFT-type if, in addition,  the adjoint module 
$V$ is isomorphic to its dual $V'$.
For a state $u \in V$, let
\begin{eqnarray*}
Y(u,z)=\sum_{n\in\Z}u_nz^{-n-1}
\end{eqnarray*}
denote the vertex operator corresponding to $u$.
It is  well-known that for such $V$,
the bracket $[u,v] =u_0v$ and the pairing $\langle u, v \rangle \1 = u_1v$
equip the weight $1$ subspace $V_1$ with the structure of  Lie 
algebra together with
a symmetric, invariant, bilinear form.  It is further proved
in \cite{DM2}, under some rationality assumptions,  that $V_1$ is
   \emph{reductive}.
   It is also well-known that the component operators of
$u_n$ ($u\in V_1, n\in\Z$) generate an \emph{affine Lie algebra}
of operators acting on $V.$  In this way, $V$ becomes a module over 
the affine Lie algebra.
Suppose that $\g$ is a simple Lie subalgebra of $V_1$.
The corresponding component operators $u_n$  generate
a copy of an \emph{affine Kac-Moody algebra} $\hat \g$, and the vertex operator
subalgebra $U \subseteq V$ generated by $\g$ is a highest weight module for
$\hat\g.$  The main result of the present paper (Theorem \ref{t3.1}) 
is the following:
\begin{eqnarray}\label{statement}
&& \mbox{\emph{Suppose that $V$ is of strong CFT-type and 
$C_2$-cofinite. Then the}} \nonumber \\
   &&\mbox{\emph{restriction of $\langle \  , \rangle$ to $\g$
   is nondegenerate, $V$ is an integrable} } \\
   &&\mbox{\emph{$\hat\g$-module,  and there is a positive integer $k$ 
such that $U \cong L(k, 0).$}} \nonumber
\end{eqnarray}

The integrality of $k$ and the structure of the integrable 
$\hat{\g}$-modules places  restraints on the structure of $V$. 
(\ref{statement})  illustrates the remarkable fact
   that  vertex operator algebras satisfying suitable finiteness constraints
   ($C_2$-cofiniteness and/or rationality), have strong 
\emph{arithmetic} properties.

    The examples of pathologies that can be exhibited by vertex operators,
   discussed in \cite{DM4}, show that we cannot weaken the assumption
   of strong CFT-type in our main result. The concept of 
$C_2$-cofiniteness, which we discuss in more detail in Section 2,
    was introduced by Zhu  \cite{Z} and has played a fundamental r\^{ole}
    in vertex operator algebra theory ever since.
  Interestingly, one of the basic features of $C_2$-cofinite vertex 
operator algebras,
    namely that
  they possess a PBW-type spanning set
with the nonrepeat condition  (\cite{GN}, \cite{B}, \cite{ABD}), was
    not discovered until several years after the publication of Zhu's paper.
  Many important results in the field seem to require
both $C_2$-cofiniteness and rationality, which together are equivalent
to regularity (\cite{Li4},  \cite{ABD}).  Our main result
    can be used to study this situation, at least for vertex operator algebras
    satisfying $V_1 \neq 0$. As an immediate corollary of the main 
result, for example,
we show that a vertex operator algebra associated to a highest weight
module for an affine Kac-Moody algebra is rational if,  and only if,
it is $C_2$-cofinite. So rationality and $C_2$-cofiniteness
are \emph{equivalent} for such vertex operator algebras. 
Understanding the extent
  of this equivalence is currently a matter of considerable interest. 
This and some other consequences are covered in Section 5 following 
the proof of the main Theorem in Section 4.

It is natural to also consider a variation of (\ref{statement}) in which
  $\g$ is replaced by an \emph{abelian} Lie subalgebra, and we discuss
  this case in Section 6. Roughly speaking,  we expect
  that $U$ can be replaced by a \emph{lattice} vertex operator algebra.
  We make an explicit Conjecture in Section 6, and show
  (Theorem \ref{tt1}) how this can be proved with an additional assumption.

\section{Background}\label{section2}
  \setcounter{equation}{0}

We start with  a simple vertex operator algebra  $V$ of
  CFT-type. Recall  the following definition \cite{Z}:   $V$ is  $C_2$-cofinite
if $V/C_2(V)$ is finite dimensional,  where $C_2(V)$
is the linear span of the states $a_{-2}b$ for $a,b\in V.$

  \medskip
  Fix homogeneous states $x^{a} \in V$, ($a$
  ranging over an index set $A$),  such that the cosets
  $x^a+C_2(V), a \in A,$ span a complement to $\1+C_2(V)$ in $V/C_2(V).$
  Set $X = \{x^a\}_{a \in A}$.  Note that all states in $X$ have 
weight at least $1$.

\begin{th}\label{tgn}  \cite{GN} $V$  is spanned by $\1$ together with
elements of the form
$$x^{1}_{-n_1} x^{2}_{-n_2} \cdots x^{k}_{-n_k} \mbox{\bf 1}$$
where $n_1>n_2> \cdots >n_k > 0,$ $x^{i} \in X$.
\end{th}

\begin{rem}\label{remsec2}   It is a consequence of the  Theorem 
that $V$ is  \emph{finitely generated} if it is $C_2$-cofinite.
\end{rem}

Next we review some of the theory of affine Kac-Moody algebras and 
their integrable
modules \cite{K}. Let $\g$ be a finite-dimensional simple
Lie algebra, $\h \subseteq \g$ a Cartan subalgebra,
  $\Phi$ the associated root system with simple roots $\Delta$, and 
$(\cdot ,\cdot)$  the nondegenerate, symmetric, invariant
bilinear form on $\g$ normalized so that the longest positive root 
$\theta \in \Phi$ satisfies
  $(\theta,\theta)=2.$   The corresponding affine Kac-Moody
algebra $\hat\g$ is defined to be
$$\hat\g=\g\otimes\C[t,t^{-1}]\oplus \C K,$$
where $K$ is central and
\begin{eqnarray} \label{LAbrackets}
  [a(m),b(n)] = [a,b](m+n)+m\delta_{m+n,0}(a,b)K,  \ (a, b \in \g, \ 
m, n \in \Z).
  \end{eqnarray}
Here,  $a(m)=a\otimes t^m.$

  \medskip
    Highest weight irreducible
  $\hat{\g}$-modules are parametrized by
  linear forms on  $\h \oplus{\bf C}K.$ These in turn correspond to pairs
  $(\lambda, k)$ where $\lambda \in \h^*$ is a weight and $K \mapsto k.$
  The corresponding $\hat{\g}$-module is denoted by
  $L(k, \lambda).$ Note that $K$ acts on $L(k, \lambda)$ as multiplication
  by $k$.  Given any $\hat{g}$-module $M$, we say that $M$
  has \emph{level $k$} if $K$ acts on $M$ as multiplication by $k$. 
Thus $L(k, \lambda)$ has level $k$.

  \medskip
  $L(k, \lambda)$ is a quotient of the Verma module
$V(k, \lambda) = U(\hat \g)\otimes_{U(\g\otimes \C[t]\oplus \C K)}\C$, ie.
the induced module Ind$_{U(\g\otimes \C[t]\oplus \C K)}^{U(\hat \g)} \C$
  for which $\g\otimes \C[t]$  acts on $\C$ via $\lambda$
  and $K$ acts as multiplication by $k.$  One knows that 
$L(k,\lambda)$ \emph{integrable}
if, and only if, $k$ is a nonnegative integer and $\lambda$ is a dominant
integral weight.

  \medskip
  For any scalar $k$ \emph{not equal} to minus the dual Coxeter 
number, $L(k,0)$ carries the structure of a
vertex operator algebra
  of strong CFT-type. It is generated by the weight $1$
  subspace $L(k,0)_1,$ which is isomorphic (as Lie algebra)
to $\g.$ The irreducible modules over the vertex operator algebra
  $L(k, 0)$ are precisely the integrable highest weight modules
$L(k,\lambda)$ satisfying $(\lambda,\theta)\leq k.$
We refer the reader to \cite{DL}, \cite{FZ}, \cite{Li3}, \cite{MP} 
for further background
  concerning these results.

  \medskip
  Finally, a word on notation. We will often be considering vertex 
operator subalgebras
  $U \subseteq V$. Generally, this does
  \emph{not} imply that $U$ and $V$ have the same Virasoro element. If they
  \emph{do}  have the same Virasoro element, we say that $U$ is a
  \emph{conformal} subalgebra of $V$.

  \section{The main Theorem}
  \setcounter{equation}{0}

We now  assume that $V$ is a vertex operator algebra that is
   \emph{simple, $C_2$-cofinite, of strong CFT-type}.
  Then $V_1$ carries the structure of a Lie algebra
with bracket $[a,b]=a_0b,  a,b\in V_1$. Because of the hypothesis
   that $V$ is of strong CFT-type, a theorem of  Li \cite{Li1} 
guarantees the existence
   and uniqueness of a symmetric, invariant bilinear form $\langle\ , 
\ \rangle$ on $V$ satisfying
  $\langle \1, \1 \rangle = -1.$  Furthermore, $\langle \ , \ 
\rangle$ is nondegenerate thanks
  to the simplicity of $V$. The restriction of this form to $V_1$ is 
thus also nondegenerate,
  and satisfies
  $a_1b =  \langle a, b \rangle \1, \  a, b \in V_1.$

  \medskip
For states $u, v \in V_1$
  and integers $m, n$ we have
\begin{eqnarray}\label{a1}
[u_m, v_n]  =  (u_0v)_{m + n}  +  m u_1v \delta_{m+n, 0}.
\end{eqnarray}
Comparison with (\ref{LAbrackets}) shows that
  if $\g \subseteq V_1$ is a (simple) Lie subalgebra
  with the property that the restriction of $\<u,v\>$ to $\g$
  is nondegenerate, then the map
  \begin{eqnarray*}
  u(m)\mapsto u_m, \ \ u\in \g, m\in\Z,
\end{eqnarray*}
is a \emph{representation of the affine Kac-Moody algebra $\hat \g$ of level
$k$},  where
\begin{eqnarray}\label{leveldef}
  \<u,v\>=k(u, v), \ \ u,v\in \g.
\end{eqnarray}
  In this way,  $V$ acquires the structure of a $\hat{g}$-module of 
level $k$, although some care should be taken with this statement. If 
$U \subseteq V$ is the vertex operator subalgebra generated by 
$\hat{g}$,
  then $U$ and $V$ generally have \emph{distinct} Virasoro elements. 
As a result,
  the $L(0)$-grading on $V$ is generally \emph{not} the same as the 
grading which arises from
  the $\hat{g}$-action.
   We can now state the main result of the present paper:

\begin{th}\label{t3.1}  Let  $V$ be a  simple vertex operator algebra 
$V$ which is
$C_2$-cofinite of strong CFT type, with $\g \subseteq V_1$
   a simple Lie subalgebra, $k$
  the level of $V$ as $\hat{g}$-module, and $U$  the vertex operator 
subalgebra of $V$
  generated by $\g.$ Then the following hold:
  \begin{eqnarray*}
  &&(a) \ \mbox{The restriction of $\langle \ , \rangle$ to $\g$ is 
nondegenerate}, \\
&&(b) \ U \cong L(k, 0), \\
&&(c) \ k \ \mbox{is a positive integer}, \\
&&(d) \ V \ \mbox{is an integrable $\hat{g}$-module}.
\end{eqnarray*}
\end{th}

\pf We first prove the Theorem in case $\g \cong sl(2, \C)$, and then 
extend it to
the general case. Until further notice, then, assume that $\g \cong sl(2, \C)$
with standard basis $\alpha,x_\alpha,x_{-\alpha}.$
   Then $(\alpha, \alpha) = 2$, so that by (\ref{leveldef}) we have
  \begin{eqnarray}\label{kvalue}
k =  \<\alpha,\alpha\>/2.
\end{eqnarray}

  Note that $U$ itself is
a highest weight module for $\hat
   \g,$ indeed $U$ is a quotient of the Verma module
$V(k,0).$ Thus
   $U$ is integrable if, and only if,
$U \cong L(k,0)$ and $k$ is a
   nonnegative integer.
By Theorem 10.7 in \cite{K}, this is in turn equivalent to
  the condition
  \begin{eqnarray}\label{nilpotentcond}
   (x_\alpha)_{-1}^l\1=0, \ \ \mbox{for some} \ l \geq 0.
   \end{eqnarray}
   By Proposition 13.6 in \cite{DL},
  $V$
   is integrable (as $\hat\g$-module) if, and only
if, $U$ is
   integrable. Thus in order to prove parts (b) and (d) of the Theorem,
   it suffices to show that  (\ref{nilpotentcond})
   holds.  Once this is achieved, it further follows (since $U$ 
contains $\g$ and hence
   cannot be $1$-dimensional) that $k$ is in fact positive and part 
(c) holds. Then (\ref{kvalue}) tells
us that (a) also holds. The upshot is that we only need prove 
(\ref{nilpotentcond}).

\medskip
Now each homogeneous subspace of $V$ is a
   $\g$-module, so that $V$ is a completely reducible $\g$-module. Of
   course,
$x\in\g$ acts as $x_0.$ A nonzero element $v\in V$ is
   called a \emph{weight vector
(for $\g$) of $\g$-weight $\lambda \
   (\lambda\in \C\alpha$)} if
$\alpha_0v=(\alpha,\lambda)v.$ From the
   representation theory
for $sl(2,\C)$ we know that $\lambda$ is a
   weight if, and only
if, $\lambda\in\frac{1}{2}\Z\alpha=P.$ Define
   an order on $P$
according to the usual order on $\frac{1}{2}\Z.$

\medskip

Let $X$ be as in Section \ref{section2}.  Since $C_2(V)$ is a 
$\g$-submodule of $V$
  we may, and shall, assume that  each  $x \in X$ is a weight
vector. Since $X$ is a finite set there is  a nonnegative
   element $\lambda_0=m\alpha \in P$
  such that the weight
of each $x\in X$ is bounded above by $\lambda_0.$

   \medskip
According to Theorem \ref{tgn}, $V$ is spanned by $\1$ together with 
elements of the form
\begin{eqnarray}
   u=x^{1}_{-n_1} x^{2}_{-n_2} \cdots x^{s}_{-n_s} \mbox{\bf 1} \label{ustate}
   \end{eqnarray}
where $n_1>n_2> \cdots >n_s > 0$ and $x^{i} \in X.$
   Since the
  $L(0)$-weight of each $x\in X$ is at least one then $L(0)$-weight of 
$u$ is at least $s(s+1)/2.$ It follows that
   for any integer $t \geq 0$,
     \begin{eqnarray}\label{Vncontainment}
\bigoplus_{n \leq t(t+1)/2} V_n \subseteq \langle x^{1}_{-n_1} 
x^{2}_{-n_2} \cdots x^{s}_{-n_s} \mbox{\bf 1}
\ | \ n_1>n_2> \cdots >n_s> 0, x^i \in X, s \leq t \rangle.
\end{eqnarray}
   On the other hand,
since each $x \in X$ is a weight vector then so too is
  the state $u$ in (\ref{ustate}), and the $\g$-weight of $u$ is at 
most $s\lambda_0.$
  From (\ref{Vncontainment}) it follows that
  \begin{eqnarray}
  \mbox{If $n \leq t(t+1)/2$, a  $\g$-weight vector in $V_n$} \\  \nonumber
  \mbox{has $\g$-weight  no greater than $t\lambda_0 = tm\alpha$}. 
\label{gweight}
  \end{eqnarray}

Now choose an integer $l$ such that $(l+1)/2 > m,$
  and set $u_0 = (x_{\alpha})_{-1}^{l(l+1)/2} \1.$ We will show that
   $u_0 = 0$, so assume that this is false.
   Then the $L(0)$-weight of $u_0$ is
   $l(l+1)/2$, so that  by (\ref{gweight}) its $\g$-weight
    is no greater than $lm \alpha.$ On the other hand, direct 
calculation shows that
    the $\g$-weight of $u_0$ is $(l(l+1)/2) \alpha > lm \alpha.$
  This contradiction shows that indeed $u_0 = 0$
    and, as explained before,  completes the proof of the Theorem in case
    $\g \cong sl(2, \C).$

\medskip

  Turning to the general case, we have the root space decomposition
\begin{eqnarray*}
    \g=\h \oplus \oplus_{\alpha \in \Phi} \g_{\alpha}.
    \end{eqnarray*}
Let $\g_{\alpha}=\C x_{\alpha}.$ The discussion in the second paragraph of
    the Proof applies to general $\g$, and is modified only in that in order to
    conclude
  $U \cong L(k,0)$, we must establish (\ref{nilpotentcond}) for
  \emph{all}  $\alpha\in\Phi.$
But for any such $\alpha$ we can find $h_{\alpha} \in\h$ such that
    $\C h_{\alpha} +\C x_{\a}+\C x_{-\alpha}$
is isomorphic to $sl(2,\C),$ and (\ref{nilpotentcond}) follows by the 
case already proved.
  This completes the proof of the Theorem. \ \ \ \ \ \
\qed

\bigskip
\section{Consequences of the Main Theorem}
\setcounter{equation}{0}

There are a number of consequences of the main Theorem. We collect a 
few of them in this Section. We retain previous notation.  The first 
result concerns the relation between
   \emph{weak, admissible} and
\emph{ordinary} $V$-modules, for which we refer to \cite{DLM1}.

\begin{coro}\label{c22}
  Let the assumptions and notation be as in Theorem \ref{t3.1}.  Then 
any weak $V$-module is an integrable $\hat \g$-module.
  \end{coro}

\pf By Theorem \ref{t3.1},   $U \cong L(k, 0)$ for some positive integer $k.$
  Now it is clear that a weak $V$-module is also a weak $U$-module. Moreover, it
is proved in \cite{DLM1} (cf. \cite{ABD})
that $L(k,0)$ has the property
that any weak module is a direct sum of ordinary irreducible $L(k,0)$-modules.
We have already mentioned that the inequivalent irreducible
$L(k,0)$-modules are exactly the level $k$ integrable highest weight
modules for $\hat\g,$ so that a weak $V$-module, considered as $U$-module,
  is the direct sum
of  integrable highest weight modules of level $k$, and is thus 
itself integrable.  This completes
the proof of the Corollary. \ \ \ \
\qed

\medskip
  Next  we have the following well-known result in the literature.
\begin{coro}\label{c2} Let $\g$ be a finite-dimensional simple Lie algebra.
Suppose that $W$ is a highest weight $\hat{g}$-module
  with highest weight $0$ with the property that
$W/(\g\otimes t^{-2}\C[t^{-1}])W$ is finite-dimensional.
  Then $W \cong L(k,0)$ for some positive integer $k$.  In particular,
$W$ is integrable.
\end{coro}

\pf Because $W$ is a highest weight module of highest weight $0$
then it is a quotient of the Verma module $V(k, 0)$ for some $k$.
  In particular, $W$ is itself a vertex operator algebra. Now it is 
proved in \cite{DLM1} that $C_2(W)=(\g\otimes t^{-2}\C[t^{-1}])W,$ so 
that the finiteness hypothesis in
  the Corollary amounts to saying that $W$ is $C_2$-cofinite.
  The Corollary now follows from Theorem \ref{t3.1} by taking $U = V = 
W.$ \ \ \ \ \ \
\qed

\medskip
If we also assume that $V$ is
\emph{rational},  then the main Theorem has particularly striking 
consequences. It is
known [ABD] that the conjunction of $C_2$-cofiniteness and 
rationality is equivalent
to so-called \emph{regularity} (loc. cit.) Thus we say that a vertex 
operator algebra
$V$ is \emph{strongly regular} if it is rational and satisfies the 
other conditions imposed on $V$ in Theorem \ref{t3.1}. We then have

\begin{coro}\label{c3} Suppose that $V$ is a simple, strongly regular 
vertex operator algebra.
  Then the Lie algebra $V_1$ is reductive, and a direct sum
  \begin{equation}\label{2.21}
V_1  =  {\frak g}_{1, k_1}  \perp {\frak g}_{2, k_2} \perp ... \oplus
{\frak g}_{n,k_n} \perp  \h
\end{equation}
  which is orthogonal with respect to $\langle \ , \rangle,$ where 
$\h$ is abelian
  and each $\g_{i, k_i}$ is a simple Lie algebra whose level $k_i$ is 
a positive integer.
  The vertex operator subalgebra of $V$ generated by the semisimple 
part of $V_1$
  is isomorphic to
  \begin{eqnarray*}
L(k_1, 0) \otimes ... \otimes L(k_n, 0).
\end{eqnarray*}
\end{coro}

  \pf That $V_1$ is a reductive Lie algebra and an orthogonal direct 
sum under the stated conditions on $V,$ is Theorem 1.1 of 
\cite{DM2}. However, there is no discussion of the possible levels 
$k_i$ of the simple summands $\g_i$ there, a lacuna which is now
  filled by Theorem \ref{t3.1}. Indeed, Theorem \ref{t3.1} applies
  with $\g_i$ in place of $\g$.
  The various assertions of the Corollary follow immediately.
  \ \ \ \ \ \  \qed

  \medskip

\begin{rem} The integrality of the levels $k_i$ was already
  announced at the end of
  \cite{DM3}. As explained there, it can be used to significantly 
restrict the possibilities for the structure of $V_1$ under the 
assumption that $V$ is holomorphic of central charge $24.$
  Indeed, the results of Corollary \ref{c3} provide justification for 
some of the main assumptions
  in the work of Schellekens \cite{Sc} on the same subject, as well as 
other papers in the physics
  literature.
\end{rem}

\section{Lattice vertex operator subalgebras and abelian subalgebras of $V_1$}
\setcounter{equation}{0}

   In our discussion so far
   we have concentrated
  on the simple Lie subalgebras of $V_1$. It is natural to ask
  about the  case of (nondegenerate) abelian Lie subalgebras $\h$
  too. In the situation of Corollary \ref{c3}, for example, this amounts
  considering the abelian sumand of $V_1$.  Of course, the vertex 
operator subalgebra
  generated by $\h$ is a Heisenberg algebra, but one can expect more.

   \medskip
  \noindent
  {\bf Conjecture} There is a positive-definite even  lattice $L 
\subseteq \h$ so that
  $V$ contains a vertex operator algebra $U$ with the following properties:
   $U$ is (isomorphic to) the lattice theory $V_L$, and $\h \subseteq U.$

   \medskip
  The Conjecture is related to the methods and results of \cite{DM2}.
  Here we make a start  on  the Conjecture, using some related ideas.
   We  continue to assume that $V$ is strongly regular, although the 
first few Lemmas below do not require rationality. We take
$\h \subseteq V_1$ to be an abelian subalgebra such that the restriction
  of $\langle \ , \rangle$ to $\h$ is nondegenerate. We let $l$ be the 
rank of $\h$.
    Note that  by
\cite{DM2}, the regularity of $V$ implies that
  for any state $u \in \h$ and any irreducible $V$-module $W$, the
  zeroth component operator $u_0$ (acting on $W$) is \emph{semisimple}.

  \medskip

  Until further notice we fix $W$ to be an irreducible $V$-module. The 
component operators of the vertex
operators $Y_W(u,z)$ acting on $W$ ($u \in \h$), close on an affine Lie algebra
$\hat\h$. That is, they satisfy the relations
$$[u_m, v_n]  =  m\delta_{m,-n}\<u , v\>.$$

  We also have
\begin{eqnarray*}
M(1)=\C[u_{-m}|u\in {\frak h},m>0],
\end{eqnarray*}
  the (rank $l$) Heisenberg subalgebra of $V$ generated by $\h$. The 
tensor decomposition
\begin{equation}\label{4.1}
W =  M(1)\otimes_{\C} \Omega_W
\end{equation}
is well-known (cf. [FLM, Theorem 1.7.3]). Here,
  \begin{eqnarray*}
  \Omega_W = \{ w \in W \ | \ u_nw = 0, u \in \h, n > 0 \}
  \end{eqnarray*}
is the so-called {\em vacuum space}.
   For $\alpha \in \h,$ set
   \begin{eqnarray*}
   W(\alpha) &=& \{ w\in W \ | \ u_0w = \<\alpha,u\>w, u\in\h \}, \\
   \Omega_W(\alpha) &=& W(\alpha)\cap \Omega_W.
   \end{eqnarray*}
   Setting
  \begin{eqnarray*}
P_W = \{ \alpha \in \h \ | \ W(\alpha) \not = 0 \},
\end{eqnarray*}
  we have the following equalities:
    \begin{eqnarray*}
    W(\alpha) &=& M(1)\otimes \Omega_W(\alpha), \\
    \Omega_W &=& \oplus_{\alpha\in P_W} \ \Omega_W(\alpha), \\
    W &=& \oplus_{\alpha\in P_W} \ W(\alpha).
    \end{eqnarray*}
    (The latter two equalities hold thanks to the assumed 
semisimplicity of the operators
     $u_0$.)
     We set
     \begin{eqnarray*}
P = P_V.
\end{eqnarray*}
   From the proof of Theorem 1.2 in \cite{DM2}, $P$ is an 
\emph{additive subgroup of $\h$}.

\begin{lem} \label{lem6.1} The following hold:
     \begin{eqnarray*}
     &&(a) \ \mbox{There is $\lambda_W \in \h$ such that $P_W=P+\lambda_W$}, \\
&&(b) \  \mbox{
$W=\oplus_{\alpha\in P}W(\alpha)=\oplus_{\alpha\in P}M(1)\otimes 
\Omega_W(\alpha).$}
     \end{eqnarray*}
\end{lem}

\pf Choose any (nonzero) $w \in W$ such that $w$ is an $\h$-weight 
vector corresponding to
     $\lambda_W \in \h$. Thus, $h_0w=\langle \lambda_W,h \rangle w, \ 
h\in \h.$ Now let
     $\alpha \in P$, with $0 \not= v \in V(\alpha).$ For $h\in \h$
  we  have
$$[h_0,Y_W(v,z)]=Y_W(h_0v,z)=\langle \alpha,h \rangle Y_W(v,z).$$
     It follows that
     \begin{eqnarray*}
h_0(v_nw) =  (\langle \lambda_W, h \rangle +  \langle \alpha, h \rangle )v_nw,
\end{eqnarray*}
  and since  $Y_W(v,z)w\ne 0$ (Proposition 11.9 of \cite{DL})
  then we can conclude that $v_nw\in W(\alpha+\lambda_W)$ for $n\in\Z.$

     However, since $W$ is irreducible then we know (Corollary 4.2 of 
\cite{DM1}) that
$W$ is spanned by $v_nw$ for $v\in V$ and $n\in \Z.$ It therefore follows that
     every element of $P_W$ has the form $\alpha + \lambda_W$ for some
choice of $\alpha$.
  Both parts of the Lemma follow from this. \ \ \ \ \ \   \qed

\medskip
The next Lemma follows from  the proof of Theorem 4.4 of \cite{DM1}.
\begin{lem}\label{lem6.2}   The following hold:
\begin{eqnarray*}
&&(a) \  V(0) = M(1)\otimes \Omega_V(0) \ \mbox{is a simple, 
conformal vertex operator
subalgebra} \\
&& \hspace{0.6cm} \mbox {of $V$ }. \\
&&(b) \ \mbox{If} \ \alpha \in P_W \ \mbox{then} \ W(\alpha) = 
M(1)\otimes \Omega_W(\alpha)
\mbox{is an irreducible $V(0)$-module.}\\
&& \hspace{0.6cm} \mbox{Moreover} \ W(\alpha) \not\cong W(\beta) \ \mbox{if} \
    \alpha \ne \beta. \ \ \ \ \ \ \ \ \ \ \ \ \ \ \ \ \ \ \ \qed
\end{eqnarray*}
\end{lem}

Let $\lambda\in\h.$ Recall that $M(1,\lambda)=M(1)\otimes e^{\lambda}$ is the
irreducible
$M(1)$-module such that $h_0$ acts as $\langle h,\lambda \rangle$ for
$h\in\h.$  Let $\{h^1,...,h^l \}$ be an orthonomal basis of $\h$ with respect
to $\langle \ , \rangle$. Then
\begin{eqnarray}\label{M(1)Virasoro}
\omega_{\h}=\frac{1}{2}\sum_{i=1}^l (h^i_{-1})^2
\end{eqnarray}
  is the Virasoro vector for $M(1).$
The following  is immediate:
\begin{coro}\label{cor6.3}
  $\Omega_V(0)$ is a simple vertex operator subalgebra of
$V$ with Virasoro element $\omega-\omega_{\h}.$ Moreover, as a 
$V(0)$-module, $V(\alpha)=M(1,\alpha)\otimes \Omega_{V}(\alpha)$
and each $\Omega_V(\alpha)$ is an irreducible $\Omega_V(0)$-module. \ 
\ \ \ \ \ \ \ \ \ \ \ \ \qed
\end{coro}

   Let us now further assume that $V$ is strongly regular, as 
previously defined.
    From the proof
of Theorem 1 in \cite{DM2} we know that $h_0=0$ for
$h\in\h$ if and only if $h=0.$ Thus $P$ spans $\h.$
\begin{lem}\label{lp} $P$ is a free abelian group of  rank at least $l$.
\end{lem}

\pf Since $V$ is $C_2$-cofinite, it is finitely generated
   (cf. Remark \ref{remsec2}). This implies that $P$ is a free abelian 
group of finite
rank. Since $P$ spans $\h,$ the rank of $P$ is at least $l.$ \ \ \ \ \ \ \
\qed.

\medskip
We now  have the following contribution to the Conjecture:
\begin{th}\label{tt1} Let $V$ be a strongly regular simple vertex 
operator algebra
  and let $\h \subseteq V_1$ be a subalgebra of rank $l$ contained in 
a Cartan subalgebra
and such that the restriction of $\langle \ , \rangle$ to
  $\h$ is nondegenerate. Assume that
the rank of $P$ is \emph{exactly} $l$. Then there are a 
positive-definite even lattice $L \subseteq P$
and a vertex operator subalgebra $U \subseteq V$ with $|P: L| < \infty$ and
  $\h \subseteq U \cong V_L.$
\end{th}

\pf Since $P$ has rank $l$ by assumption, and since the restriction of
  $\langle \ , \rangle$ to $\h$ is nondegenerate, it follows that the
  $\Z$-dual lattice $P^{\circ}$ of $P$  also
  has rank $l.$

   Now let $W$ be a $V$-module and $\alpha\in P^{\circ}.$
We define a new action of $V$ on $W$ by
$$Y_{W}^{\alpha}(v,z)=Y_W(\Delta(\a,z)v,z)$$
for $v\in V,$ where
$$\Delta(\a,z)=z^{\a(0)}\exp\left(\sum_{k=1}^{\infty}\frac{\a(k)}{-k}
(-z)^{-k}\right)$$
is Li's $\Delta$-field \cite{Li2}. Then $W^{\alpha}=(W,Y_W^{\alpha})$ is again
   a $V$-module, moreover $W$ is irreducible if, and only if, 
$W^{\alpha}$ is irreducible (loc. cit.) Furthermore, it is easy to 
see that
$(W^{\a})^{\b} \cong W^{\a+\b}$ for $\alpha,\beta\in P^{\circ}.$

Let $\irr(V)$ be the set consisting of the isomorphism classes
   of inequivalent irreducible $V$-modules. The $C_2$-condition and the
   rationality assumption both imply that $\irr(V)$ is a \emph{finite} set
   (\cite{DLM2}, \cite{Z}). For an irreducible $V$-module
   $W$, let $[W]$ denote the element in $\irr(V)$ determined by $W$. Simplicity
   of $V$ says that
$[V] \in \irr(V).$
It follows from what we have said that there is a natural group action
   $P^{\circ} \times \irr(V) \rightarrow \irr(V)$, namely
   \begin{eqnarray*}
(\alpha,  [W]) \mapsto  [W^{\alpha}].
\end{eqnarray*}
We define $L$ via
\begin{eqnarray*}
L = \ \mbox{stab}_{P^{0}}([V]) = \{ \alpha \in P^0 \ | \ V^{\alpha} \cong V \}.
\end{eqnarray*}
The finiteness of $\irr(V)$ implies that $|P^o:L| < \infty$,  so
that $L$ is a free abelian group of rank $l.$

Suppose that $\alpha\in L.$ Then there is an isomorphism of $V$-modules
$V \stackrel{\cong}{\rightarrow} V^{\alpha}$. This amounts to the existence
of a linear isomorphism $\phi_{\alpha}:V\to V$ satisfying
\begin{eqnarray*}
  Y(\Delta(\a,z)v,z)=\phi_{\a}Y(v,z)\phi_{\a}^{-1}, \ v \in V.
\end{eqnarray*}
  In particular, $Y(v,z)= \phi_{\a}Y(v,z)\phi_{\a}^{-1}$
for $v\in \Omega_V(0).$ That is,
$\phi_{\alpha}$ is an $\Omega_V(0)$-module isomorphism.

Note that $\Delta(\a,z)h(-1)=h(-1)+(\alpha,h)z^{-1}$ for
$h\in \h.$ This shows that $\phi_{\alpha}(V(\beta))=V(\beta+\alpha)$
for $\beta\in P.$ In particular, $V(\beta)$ and $V(\beta+\alpha)$ are
isomorphic $\Omega_V(0)$-modules. As a result, there is an isomorphism
  of $\Omega_V(0)$-modules
  \begin{eqnarray}\label{Omegaisom}
\Omega_V(\beta) \cong \Omega_V(\beta+\alpha), \ \mbox{all} \ \alpha \in L.
  \end{eqnarray}
   Taking $\beta = 0$, it follows that
  $\Omega_V(\alpha) = \Omega_V(0)$ for all $\alpha \in L$, and therefore
  $L \subseteq P$.  Set $|P:L| = p$ and choose a set of coset representatives
  $\mu_1, ..., \mu_p$ of $L$ in $P$ with $\mu_1 = 0.$
  Then (\ref{Omegaisom}) together with Lemma \ref{lem6.1}(b) and Corollary
  \ref{cor6.3} shows that, when considered as a module
  over the vertex operator subalgebra $M(1)\otimes \Omega_V(0)$,
  we have a decomposition
$$V \cong \sum_{i = 1}^{p} \sum_{\alpha\in L}M(1,\a+\mu_i)\otimes 
\Omega_V(\mu_i)$$
of $V$ into irreducible $M(1)\otimes \Omega_V(0)$-submodules.

We are now ready to show that $L$ is a positive-definite even lattice.
  Let $\alpha \in L.$ Because
  of Lemma 6.2(a), the last display shows that the $L(0)$-weight
  of a highest weight vector $e^{\alpha}$ in $M(1, \alpha)$ is the \emph{same}
  as its $L_{{\h}}(0)$-weight (cf. (\ref{M(1)Virasoro})). Therefore the
weight in question is (using (\ref{M(1)Virasoro}))  $\langle \a,\a 
\rangle/2$, and is necessarily a nonnegative integer. Indeed it is 
positive unless $\alpha = 0$ since $V$ is of CFT-type.
  Thus $\langle \alpha, \alpha \rangle$ is a positive even integer 
whenever $\alpha \not= 0$,
  and $L$ is indeed a positive-definite, even lattice.

Invoking the fusion rules for $M(1)$-modules,
  $\oplus_{\alpha\in L}V(\alpha)=(\oplus_{\alpha\in 
L}M(1,\alpha))\otimes \Omega_V(0)$ is a vertex operator subalgebra of 
$V$ with the same
Virasoro element and similarly $\oplus_{\alpha\in L}M(1,\alpha)$ is a vertex
operator subalgebra with Virasoro element $\omega_{\h}.$ By Corollary
5.4 of \cite{DM2}, it follows that $\oplus_{\alpha\in L}M(1,\alpha)$ 
is isomorphic to $V_L.$
This completes the proof of the Theorem. \ \ \ \ \ \ \
\qed

  \medskip
\begin{rem}\label{tt2} Theorem \ref{tt1} shows that the Conjecture is 
true if, and
  only if,
  $P$ has rank $l$. In particular, if $P$ is \emph{rational} then by
  Lemma \ref{lp} it necessarily has rank exactly $l$, so that the 
Conjecture holds
  if $P$ is rational. This result was first proved by Li
  \cite{Li5}.
In fact, Theorem \ref{tt1} shows that the Conjecture, the rationality of $P$,
  and
  the rank condition are \emph{equivalent}.
  \end{rem}

The next corollary to Theorem \ref{tt1} is well known (cf. [DL], [LP]).
\begin{coro} Let $\g$ be a finite-dimensional simple Lie algebra and  $k$
  a positive integer. The simple vertex operator algebra
$L(k,0)$ associated to $\g$ contains a lattice vertex operator 
subalgebra $V_{\sqrt{k}Q_l},$ where
$Q_l$ is the lattice generated by the long roots normalized so that 
the squared length
   is 2.
\end{coro}

\medskip
  Our finally result is
  \begin{coro} Suppose that $V$ is a strongly regular vertex operator 
algebra, and
  let $G \subseteq V$ be the vertex operator subalgebra generated by $V_1$. Then
  the following hold:
  \begin{eqnarray*}
&&(a) \ \mbox{$V = G$ if, and only if, $V \cong L(k_1, 0) \otimes ... 
\otimes L(k_n, 0)$  for positive integers $k_1, ..., k_n$}; \\
&&(b) \ \mbox{$G$ is a \emph{conformal} vertex operator subalgebra
  if, and only if, $G$ is contained in a} \\
&& \ \ \ \ \mbox{rational conformal vertex operator subalgebra
$H \cong L(k_1, 0) \otimes ... \otimes L(k_n, 0) \otimes V_L$} \\
&& \ \ \   \mbox { where $k_1, ..., k_n$ are positive integers
and $L$ is a
  positive-definite even lattice.}
\end{eqnarray*}
\end{coro}

\pf Note from Corollary \ref{c3} that
\begin{eqnarray}\label{Gdecomp}
G \cong L(k_1, 0) \otimes ... \otimes L(k_n, 0) \otimes M(1),
\end{eqnarray}
where $M(1)$ is the Heisenberg vertex operator subalgebra generated by
the abelian summand $\h$ in (\ref{2.21}).  Part (a) follows immediately.
  As for part (b), suppose that $G$ is a conformal vertex operator
subalgebra of $V$, so that $\omega$ is the Virasoro vector of $G$.  Since
the $k_i$ are positive integers then the weight spaces of the modules over
$L(k_1, 0) \otimes ... \otimes L(k_n, 0)$ have rational weights. It 
then follows
that $P$ is a rational lattice, and by Remark \ref{tt2} we may apply Theorem
\ref{tt1} to see that there is a tower of conformal vertex operator subalgebras
$G \subseteq H \subseteq V$ with $H$ as in part (b). Conversely, if 
$G$ is contained
in a conformal vertex operator subalgebra $H \cong   L(k_1, 0) 
\otimes ... \otimes L(k_n, 0) \otimes V_L$, then $G$ is necessarily a 
conformal vertex operator subalgebra of $H$, and hence also of $V$. \ 
\ \ \ \
\qed

\end{document}